\input amssym.def
\input amssym
\def\beginlist#1{\begingroup\parindent=#1}
\def\endlist{\smallskip\noindent\endgroup}
\def\litem{\smallskip\noindent
\hangindent=\parindent\ltexindent}
\def\ltexindent#1{\hbox to \hangindent{#1\hss}\ignorespaces}
\magnification=\magstep1
\def\nbigskip{\bigskip\noindent}
\def\bbigskip{\bigskip \bigskip  }
\def\nbbigskip{\bbigskip \noindent }
\def\nbbbigskip{\bbigskip \bigskip \noindent }
\def\nmedskip{\medskip\noindent}
\def\mal{\mathbin{\! \cdot \!}}
\def\qq{/\kern-.185em /}
\def\C{{\Bbb C}}
\def\red{\mathop{\hbox{\rm red}}}
\def\Hom{\mathop{\hbox{\rm Hom}}}
\def\Lin{\mathop{\hbox{\rm Lin}}}
\def\Im{\mathop{\hbox{\rm Im}}}
\def\supp{\mathop{\hbox{\rm Supp}}}
\def\buildover#1#2{\buildrel#1\over#2}
\def\qed{\hfill${\scriptstyle\square}$}
\def\pir{\pi_{\rm r}}

\centerline{\bf Actions of compact groups on coherent sheaves}

\bigskip
\centerline{\sl J\"urgen Hausen and Peter Heinzner\footnote{$^1$}
{Supported by a Heisenberg Stipendium of
the Deutsche Forschungsgemeinschaft}}

\nbbbigskip\goodbreak\noindent
Let $K$ be a compact Lie group with complexification $K^\C$ and 
$X$ a reduced Stein space endowed with a continuous action of $K$
by holomorphic transformations. In this set up there is a complex 
space $X\qq K$ and a surjective $K$-invariant holomorphic map $\pi
\colon X \to X \qq K$ such that the structure sheaf of $X \qq K$ is
the sheaf of invariants $(\pi_*{\cal O}_X)^K$. Furthermore, the
$K$-action on $X$ can be complexified in the following sense: $X$ can
be realized as an open $K$-stable subset of a Stein space $X^c$ which
is endowed with a holomorphic $K^\C$-action. Moreover, the inclusion
$X\subset X^c$ is universal, i.e., every $K$-equivariant holomorphic
map $\phi$ from $X$ into a holomorphic $K^\C$-space $Z$ extends
uniquely and $K^\C$-equivariantly to a holomorphic map $\phi^c: X^c
\to Z$ (see [H]).

\smallskip
In the present paper we extend the above results to Stein spaces $X$
which are not necessarily reduced. More importantly, we show that
every continuous coherent $K$-sheaf over $X$ can be extended to a
coherent $K^\C$-sheaf over $X^c$. This is also proved in the context
of non-reduced spaces.

\smallskip
The constructions of quotients and complexifications are compatible
with the reduced structures. Thus we have $\red (X) \qq K = \red(X \qq
K)$ and $\red (X)^c = \red (X^c)$, where $\red$ associates to a
complex space its underlying reduced space. Moreover, the results are
also valid for $K$-spaces $X$ with $\red (X)$ admitting a semistable
quotient (see Section 1). In this context we show that the sheaf
$(\pi_*{\cal S})^K$ of invariant sections of a coherent $K$-sheaf
${\cal S}$ over $X$ is a coherent analytic sheaf over the quotient
space $X\qq K$ (see Section 1). Two direct applications of our results
are

\smallskip
\item{a)} For every closed complex $K$-subspace $A$ of $X$ the image
$\pi(A)$ is a closed subspace of $X\qq K$ such that $A\qq K\cong \pi(A)$.
Moreover, $A$ extends to a closed $K^\C$-subspace $A^c$ of $X^c$.

\smallskip
\item{b)} Every holomorphic $K$-vector bundle over $X$ is given by the
restriction of a holomorphic $K^\C$-vector bundle over the
complexification $X^\C$.

\smallskip
The authors would like to thank the referees for their interest
in the subject and for valuable suggestions.

\nbbigskip\goodbreak\noindent
{\bf 1. Formulation of the results}

\nbigskip
In this paper a complex space $X$ is a not necessarily reduced complex
space with countable topology. By $\red(X)$ we denote the underlying
reduced space, i.e., set-theoretically $\red(X)=X$ and ${\cal
O}_{\red(X)}= {\cal O}_X/ {\cal N}_X$ on the level of sheaves. Here
${\cal N}_X$ denotes the nilradical of the structure sheaf ${\cal
O}_X$ of $X$.

\smallskip
Let $G$ be a real Lie group and assume that $G$ acts continuously 
on $X$ by  holomorphic transformations. Then $G$ acts also on the
structure sheaf ${\cal O}:={\cal O}_X$ of $X$: A given $g\in G$
induces for every  $x\in X$ an isomorphism ${\cal O}_x\to {\cal
O}_{g \cdot x}$, $f \mapsto g\mal f$, of stalks. We call $X$ a {\it
complex $G$-space} if for all open sets $N \subset G$ and
$U_1$, $U_2 \subset X$ with $U_2 \subset \bigcap_{g \in N} g \mal
U_1$ the map
$$ \Phi\colon N \times {\cal O}(U_1) \to {\cal O}(U_2), 
\quad (g,f) \mapsto (g \mal f)\vert U_2 \leqno{(1)}$$
is continuous with respect to the canonical Frech\'et topology on 
${\cal O}(U_j)$, $j=1,2$ (see [G-R], Chapter V, $\S$6). 
If $G$ is a complex Lie group and all above maps 
$\Phi$ are holomorphic, then we say that $X$ is 
a {\it holomorphic $G$-space}.

\smallskip
Now let ${\cal S}$ be a $G$-sheaf on $X$, i.e., $G$ acts on
${\cal S}$ such that the projection ${\cal S} \to X$, $s_x \mapsto x$
is $G$-equivariant. Reformulating (1) for a coherent sheaf ${\cal S}$,
we obtain the notion of a {\it continuous} (resp. {\it holomorphic})
{\it coherent $G$-sheaf} on $X$.

\smallskip
A holomorphic map $\phi \colon X \to X'$ of two complex $G$-spaces $X$
and $X'$ is called {\it $G$-equivariant},
if it is $G$-equivariant as a map of sets and the
comorphism $\phi^0$ is compatible with the actions of $G$ on the
structure sheaves ${\cal O}_X$ and ${\cal O}_{X'}$. The map $\phi$ is
called {\it $G$-invariant} if it is equivariant with respect to the
trivial action of $G$ on $X'$.

\nbigskip
{\it Remark 1.} If $X$ is reduced, then it is a complex
$G$-space. If moreover $G$ is a complex Lie group
and the action $G\times X\to X$ is holomorphic, then $X$ is a
holomorphic $G$-space (see [K]). A holomorphic map of
reduced $G$-spaces is equivariant if and only if it is equivariant as
a map of sets.

\bigskip
Now let $K$ be a compact real Lie group and assume that $X$ is
a complex $K$-space. A $K$-invariant holomorphic map $\pi \colon X \to
X \qq K$ onto a complex space $X \qq K$ is said to be a {\it
semistable quotient} of the $K$-space $X$, if

\nmedskip
\item{i)} the structure sheaf ${\cal O}_{X \qq K}$ of ${X \qq K}$ is
the sheaf of invariants $(\pi_*{\cal O}_X)^K$, i.e., for every open
set $Q \subset X \qq K$ we have ${\cal O}_{X \qq K}(Q) = {\cal
O}_X(\pi^{-1}(Q))^K$ and

\smallskip
\item {ii)} $\pi$ is a Stein map, i.e., for every Stein open set 
$Q \subset X \qq K$ the inverse image $\pi^{-1}(Q)$ is also Stein.

\medskip
The above definition generalizes the notion introduced in [H-M-P] to
the non-reduced case. Examples of semistable quotients occur in
Geometric Invariant Theory. Our first main result is the following

\nbigskip
{\bf Quotient Theorem.} {\it A semistable quotient $\pi \colon X \to X
\qq K$ for $X$ exists if and only if it exists for $\red(X)$. If this is
the case, then $\red(\pi) \colon \red(X) \to \red(X\qq K)$ is a semistable
quotient for $\red(X)$.}

\nbigskip
By [H], Section 6.5, every reduced Stein $K$-space has a
semistable quotient and the associated quotient space is again
Stein. This, together with the above result, implies the following

\nbigskip
{\bf Corollary 1.} {\it If $X$ is a Stein $K$-space, then the
semistable quotient $\pi \colon X \to X\qq K$ exists. Moreover, $X\qq
K$ is a Stein space.}
\qed

\nbigskip
Since the reduction of a semistable quotient is a semistable
quotient of the associated reduced $K$-space, [H], 2.3, yields
that semistable quotients are universal with respect to $K$-invariant
holomorphic maps. In particular, the quotient space $X \qq K$ is
unique up to isomorphism.

\smallskip
If $X$ admits a semistable quotient $\pi \colon X \to
X\qq K$, then we denote for a $K$-sheaf ${\cal S}$ on $X$ by $(\pi_*{\cal
S})^K$ the sheaf of invariants on $X \qq K$, i.e., $(\pi_* {\cal
S})^K(Q) := {\cal S}(\pi^{-1}(Q))^K$. The following fact is a
generalization of a result of Roberts (see [R]) for holomorphic
$K^\C$-sheaves.

\nbigskip
{\bf Coherence Theorem.} {\it Assume that $X$ has a semistable quotient
$\pi \colon X \to X\qq K$ and let ${\cal S}$ be a continuous coherent 
$K$-sheaf on $X$. Then the sheaf $(\pi_*{\cal S})^K$ on $X \qq K$ is
coherent.}

\nbigskip
Let $K^\C$ denote the complexification of $K$. A {\it
complexification} $X^c$ of a complex $K$-space $X$ is a holomorphic
$K^\C$-space $X^c$ which contains $X$ as a $K$-stable open subset such
that every $K$-equivariant holomorphic map $\phi$ from $X$ into a
holomorphic $K^\C$-space $Z$ extends uniquely to a $K^\C$-equivariant
holomorphic map $\phi^c \colon X^c \to Z$. Note that a
complexification $X^c$ is unique up to  isomorphism and that  $X^c = K^\C
\mal X$.

\nbigskip
{\bf Complexification Theorem.} {\it If $X$ has a semistable quotient
$\pi\colon X\to X\qq K$, then the complexification $X^c$
exists. Moreover, the extension $\pi^c \colon X^c \to X \qq K$
is a semistable quotient for the $K^\C$-space $X^c$ and $\red(X^c) =
\red(X)^c$ holds.}

\nbigskip
This result generalizes the Theorem in Section 6.6 of [H], where the
existence of complexifications is proved for reduced Stein
$K$-spaces.

\smallskip
Now, let ${\cal S}$ be a continuous coherent $K$-sheaf on $X$ and
assume that $X$ has a semistable quotient $\pi\colon X \to X \qq K$. A
holomorphic $K^\C$-sheaf ${\cal S}^c$ on $X^c$ with ${\cal S}^c\vert X
= {\cal S}$ is called a {\it $K^\C$-extension} of ${\cal S}$.

\nbigskip
{\bf Extension Theorem.} {\it Assume that $X$ has a semistable
quotient and let ${\cal S}$ be a continuous coherent $K$-sheaf ${\cal
S}$ on $X$. Then, up to $K^\C$-equivariant isomorphism, there is a
unique $K^\C$-extension of ${\cal S}$.}

\nbigskip
For a $K$-stable closed subspace $A$ of $X$ which is  defined by a
$K$-invariant sheaf ${\cal I}$ of ideals, the Coherence and
Extension Theorem imply:

\nbigskip
{\bf Corollary 2.} {\it The sheaf $(\pi_*{\cal I})^K$ of ideals 
 endows $\pi(A)$ with the structure of closed subspace of $X\qq
K$. The restriction of $\pi$ to $A$ is a semistable quotient for
$A$. Moreover, $A$ has a complexification $A^c$ by a closed
$K^\C$-subspace of $X^c$ and $\pi(A)=\pi^c(A^c)$.}
\qed

\nbigskip
If ${\cal S}$ is a locally free continuous $K$-sheaf over $X$ with
$K^\C$-extension ${\cal S}^c$, then the complement $E$ of the set of
points in $X^c$ at which ${\cal S}^c$ is not locally free is a proper
analytic $K^\C$-stable subset of $X^c$. Since $X^c= K^\C\mal X$, this
implies that $E$ is empty. In other words we have the following

\nbigskip
{\bf Corollary 3.} {\it The extension of a locally free $K$-sheaf over
$X$ is a locally free $K^\C$-sheaf over $X^c$. In particular,
every holomorphic $K$-vector bundle over $X$ extends to a holomorphic
$K^\C$-vector bundle over $X^c$.}
\qed

\nbbigskip\goodbreak\noindent
{\bf 2. Equivariant Resolution}

\nbigskip
Let $K$ be a compact real Lie group and let $X$ be a complex
$K$-space. In this section we assume that the associated reduced
$K$-space $\red (X)$ has a semistable quotient $\pir \colon
\red (X) \to \red (X) \qq K$. Let ${\cal S}$ be a continuous coherent
$K$-sheaf on $X$. Note that we have a continuous representation of $K$
on the Frech\' et space ${\cal S}(X)$.

\smallskip
A section $s \in {\cal S}(X)$ of ${\cal S}$ is said to be {\it
$K$-finite} if the vector subspace of ${\cal S}(X)$ generated by $K
\mal s$ is of finite dimension. We denote the vector subspace 
of $K$-finite elements of ${\cal S}(X)$ by ${\cal S}(X)_{\rm fin}$. 
By a theorem of Harish-Chandra (see [HC], Lemma 5), ${\cal S}(X)_{\rm
fin}$ is dense in ${\cal S}(X)$. 

\nbigskip
{\bf Equivariant Resolution Lemma.} {\it Every $y \in \red(X) \qq K$
has an open neighborhood $Q$ such that over $U := \pir^{-1}(Q)$
there is an exact sequence ${\cal V} \to {\cal S}\vert U \to 0$ of
$K$-sheaves. Here $V$ is a finite dimensional representation space of
$K$ and ${\cal V}:={\cal O} \otimes V$ is endowed with the action
defined by $k \mal (f \otimes v) := k \mal f \otimes k \mal v$.}

\nbigskip
{\it Proof.} Since $\pir$ is a Stein map, we may assume that
$X$ is a Stein space. By Proposition 3.2 and Corollary 2 in Section
2.3 of [H], there is a unique $K$-invariant analytic set $A$ of
minimal dimension contained in $\pir^{-1}(y)$. We fix a point
$x \in A$. Then, since $X$ is a Stein space, there are sections $s_1,
\ldots, s_r \in {\cal S}(X)$ which generate ${\cal S}_x$ as an ${\cal
O}_x$-module.

\smallskip
The $K$-finite elements of ${\cal S}(X)$ are dense. Thus
we may assume that the $s_i$ generate a finite dimensional 
$K$-submodule $V$ of ${\cal S}(X)$. Hence there is an
equivariant homomorphism $\alpha\colon {\cal V} \to {\cal S}$
of $K$-sheaves which is defined by ${\bf f}\otimes s \to  {\bf f\, s}$
where ${\bf s}$ denotes the section of ${\cal S}$ defined by $s\in
V\subset {\cal S}(X)$.

\smallskip
Since $\alpha$ is a homomorphism of coherent $K$-sheaves, $B :=
\supp({\cal S} / \alpha({\cal V}))$ is a closed $K$-invariant analytic
subset of $X$. By definition, $\alpha$ is surjective over $X\setminus
B$. Furthermore, we have $\pir(A) \cap \pir(B) =
\emptyset$. Hence $Q:= \red(X)\qq K\setminus \pir(B)$ is open,
contains $y$ and $\alpha$ is surjective over $U:=\pir^{-1}(Q)$.
\qed

\nbigskip
As a consequence of the Equivariant Resolution Lemma we give here 
a proof of the Coherence Theorem in the reduced Stein case (see
also [R]). We use two well known Lemmas. The proof of the
first one can be found for example in  [R] or [Sch].

\nbigskip
{\bf Lemma 1.} {\it Let $V$ and $W$ be finite dimensional complex 
$K^\C$-modules. Then there exist $K$-equivariant polynomials 
$p_i \colon W \to V$, $i = 1, \ldots, r$ which generate 
$({\pi_W}_*{\cal V})^K$ as an ${\cal O}_{W \qq K^\C}$-module.}
\qed

\nbigskip
{\bf Lemma 2.} {\it For every exact sequence ${\cal S}_1 \to {\cal
S}_2 \to {\cal S}_3$ of continuous coherent $K$-sheaves on $X$ the
induced sequence $({\pir}_*{\cal S}_1)^K \to ({\pi_{\rm
r}}_*{\cal S}_2)^K \to ({\pir}_*{\cal S}_3)^K$ on $\red(X) \qq
K$ is also exact.}

\nbigskip
{\it Proof.} Integration over $K$ defines a projection operator
${\cal S}(X)\to {\cal S}(X)^K$. Since exactness is a local property
and $\pir \colon \red (X) \to \red(X) \qq K$ is a Stein map, the
assertion follows.
\qed

\nbigskip
{\bf Proposition 1.} {\it Let $X$ be a reduced Stein $K$-space and
denote by $\pi \colon X \to X \qq K$ the semistable quotient. Then
$(\pi_*{\cal S})^K$ is a coherent sheaf on $X \qq K$.}

\nbigskip
{\it Proof.} First consider the case ${\cal S} = {\cal V} :=
{\cal O}_{X} \otimes V$ with some $K$-module $V$. By [H], 6.6,
we may assume that $X$ is a holomorphic $K^\C$-space. Moreover, by
[H], 5.4 and 6.2, we may assume that $X$ is a closed $K$-subspace of
an open Stein $\pi_W$-saturated $K$-subspace of some finite
dimensional $K$-module $W$. So, for the case ${\cal S} = {\cal V}$,
Lemma 1 provides $p_1, \ldots, p_r \in {\cal S}(X)^K$
that generate $(\pi_*{\cal S})^K$, i.e., the associated sequence
$${\cal O}_{X \qq K}^r \buildover{\alpha}{\to} (\pi_*{\cal S})^K \to
 0.$$
is exact. The kernel of $\alpha$ is $(\pi_*{\cal R})^K$, where ${\cal
 R} \subset {\cal O}_{X}^r$ denotes the K-sheaf of relations of the
 generators $p_i$. Applying Lemma 2 to an equivariant resolution of
 ${\cal R}$ yields that also  ${\rm Ker}(\alpha)$ is finitely
 generated. Consequently $(\pi_*{\cal S})^K$ is coherent. For general
 ${\cal S}$, apply Lemma 2 to an equivariant resolution ${\cal V}' \to
 {\cal V}\to {\cal S} \to 0$. \qed

\nbbigskip\goodbreak\noindent
{\bf 3. $K^\C$-Extensions}

\nbigskip
In this section we apply the Equivariant Resolution Lemma to
investigate extensions of coherent $K$-sheaves. Assume that $X$ is a
complex $K$-space such that there is a semistable quotient $\pir
\colon \red(X) \to \red (X)\qq K$ and a complexification $X^c$ with
the following properties:

\smallskip
\item{i)} $\red (X^c)$ is the complexification of $\red (X)$ and
$\pir^c \colon \red (X^c) \to \red(X) \qq K$ is a semista\-ble
quotient.

\smallskip
\item{ii)} For every open subspace of the form $U := \pi_{\rm
r}^{-1}(Q)$ with $Q \subset \red(X) \qq K$ open, the open subspace
$U^c := K^\C \mal U = (\pir^c)^{-1}(Q)$ of $X^c$ is a
complexification of $U$.

\smallskip
Note that, by [H], Sections 3.3 and 6.6, these assumptions are valid
if $X$ is a reduced Stein $K$-space.

\nbigskip
{\bf Identity Principle.} {\it Let ${\cal T}$ be a holomorphic coherent
$K^\C$-sheaf over $X^c$. Then the restriction map $R \colon {\cal
T}(X^c)_{\rm fin} \to {\cal T}(X)_{\rm fin}$ is bijective.}

\nbigskip
{\it Proof.} First we show that $R$ is injective. So, let $s\in {\cal
T}(X^\C)$ be $K$-finite such that $s\vert X=0$. Since $u\vert
X=0$ for all $u\in V := \Lin_\C (K \mal s) = \Lin_\C (K^\C
\mal s)$, we have $(g \mal s)\vert X = 0$ for all $g\in K^\C$. Hence
$s\vert {g^{-1}\mal X} = 0$ for every $g \in K^\C$. Since $X^c =K^\C
\mal X$, it follows that $s = 0$.

\smallskip
In order to show that $R$ is surjective, we first assume that $X$ is
Stein and that there is an equivariant resolution ${\cal V}
\buildover{\alpha}{\to} {\cal T} \to 0$ with $V$ a finite-dimensional
$K$-module and ${\cal V} := {\cal O}_{X^c} \otimes V$. Every
$K$-finite $F \in {\cal V}(X)$ defines a finite-dimensional $K$-module
$W := \Lin_\C(K \mal f)$.

\smallskip
Evaluating the elements of $W$ determines a $K$-equivariant
holomorphic map $\Phi \colon X \to \Hom(W,V)$. Since $\Phi$ extends to
a $K^\C$-equivariant holomorphic map $\Phi^c \colon X^c \to V$, it
follows that $F$ extends to $X^c$ as a holomorphic map. This yields
surjectivity of the restriction ${\cal V}(X^c)_{\rm fin} \to {\cal
V}(X)_{\rm fin}$.

\smallskip
Surjectivity of $R$ is obtained as follows. Let $s \in {\cal
T}(X)_{\rm fin}$. Since $X$ was assumed to be Stein, we find an $s_1
\in {\cal V}(X)$ with $\alpha(s_1) = s$. Note that $\alpha$ maps
$\Lin_\C(K \mal s_1)$ onto the finite-dimensional vector space
$\Lin_\C(K \mal s)$. Hence we may assume that $s_1$ is $K$-finite. As
seen above, $s$ has an extension $s_1^c \in {\cal V}(X^c)_{\rm fin}$.
Then $\alpha(s_1^c) \in {\cal T}(X^c)_{\rm fin}$ is an extension of $s$.

\smallskip
Now, in the general case, let $s \in {\cal
T}(X)_{\rm fin}$. The Equivariant Resolution Lemma and the above
consideration yield a cover of $\red(X \qq K)$ by open sets
$Q_i$ such that we can extend $s$ over each $U_i = \pi_{\rm
r}^{-1}(Q_i)$ to $s_i^c \in {\cal T}(U_i^c)_{\rm fin}$. By injectivity of 
$R$, we obtain that any two such extensions $s_i^c$ and $s_j^c$ coincide over
$U_i^c \cap U_j^c$. Thus the $s_i^c$ patch together to an extension
$s^c \in {\cal T}(X^c)$ of $s$. As before we can achieve that
$s^c$ is $K$-finite.
\qed

\nbigskip
Let ${\cal S}^1$ and ${\cal S}^2$ be holomorphic coherent
$K^\C$-sheaves on $X^c$. Then ${\frak Hom}({\cal S}^1,{\cal S}^2)$ is
a holomorphic $K^\C$-sheaf. The action of $g \in K^\C$ on ${\bf F} \in
{\frak Hom}({\cal S}^1,{\cal S}^2)_x$ is given by 
$$(g \mal {\bf F})({\bf s}) := g \mal ({\bf F}(g^{-1} \mal {\bf s})). $$
where  ${\bf s} \in {\cal S}^1_{g \cdot x}$. 
The $K^\C$-invariant global sections of ${\frak Hom}({\cal S}^1,{\cal
S}^2)$ are precisely the $K^\C$-equivariant homomorphisms from ${\cal
S}^1$ to ${\cal S}^2$. So the Identity Principle yields:

\nbigskip
{\bf Homomorphism Lemma.} {\it Every homomorphism $\alpha : {\cal S}^1
\vert X \to {\cal S}^2 \vert X$ of $K$-sheaves extends uniquely to a
homomorphism $\alpha^c : {\cal S}^1 \to {\cal S}^2$ of
$K^\C$-sheaves.}
\qed

\nbigskip
{\bf Local Extension Lemma.} {\it Let ${\cal S}$ be a continuous
coherent $K$-sheaf on $X$. Then every $y \in \red(X) \qq K$ has an
open neighborhood $Q$ such that the restriction ${\cal S}_U$ of ${\cal
S}$ to $U := \pir^{-1}(Q)$ has a $K^\C$-extension ${\cal
S}_U^c$.}

\nbigskip
{\it Proof.} By the Equivariant Resolution Lemma we find an open Stein
neighborhood $Q \subset \red(X) \qq K$ of $y$  such that over $U$
there exists an equivariant resolution ${\cal V}_1\buildover \alpha
\to {\cal V}_2 \to {\cal S}_U \to 0$, where the $V_i$ are finite-dimensional
representation spaces of $K$ and ${\cal V}_i$ are the associated
$K$-sheaves on $X$.

\smallskip
Set ${\cal V}_i^c := {\cal O}_{X^c} \otimes V_i$
and endow each ${\cal V}_i^c$ with the diagonal $K^\C$-action.
By the Homomorphism Lemma, there is a unique $K^\C$-equivariant
extension $\alpha^c\colon{\cal V}^c_1\to {\cal V}^c_2$ of
$\alpha$. Identify ${\cal S}_U$ with ${\cal V}_2/ \Im(\alpha)$ and set
${\cal S}_U^c := {\cal V}_2^c/ \Im(\alpha^c)$.
\qed

\nbbigskip\goodbreak\noindent
{\bf 4. Invariant subspaces of reduced spaces}

\nbigskip
Let $K$ be a compact real Lie group and $X$ a complex
$K$-space. Assume that the associated reduced $K$-space $\red(X)$ has
a semistable quotient $\pir \colon \red(X) \to \red (X) \qq
K$. A technical ingredient for the proofs of our
results is the following

\nbigskip
{\bf Local Embedding Lemma.} {\it Every point $y \in (\red X) \qq K$ has  
an open Stein neighborhood $Q \subset (\red X) \qq K$  
such that the open subspace $U := \pir^{-1}(Q)$ of $X$ can
be realized as a closed $K$-subspace of a reduced Stein $K$-space.}

\nmedskip
{\it Proof.} We may assume that $X$ is a Stein space. As in the
reduced case (see [H], 6.2), we can find a $K$-equivariant holomorphic
map $\phi$ from $X$ into a complex finite dimensional $K$-module $V$
such that $\phi$ is an immersion along the fiber $\pir^{-1}(y)$
and $\red(\phi)$ embeds some $\pir$-saturated open neighborhood
$\red(U)$ of $x$ properly into an open $K$-stable subset $Z'$ of $V$.

\smallskip
By Siu's Theorem (see [S]), $Z'$ contains a Stein open neighborhood
$Z''$ of $\phi(U)$. Set $Z := \bigcap_{k \in K} k \mal
Z''$. Then $Z$ is an open $K$-stable Stein neighborhood of
$\phi(U)$. Now, the set $A \subset U$ consisting of all points $a \in
U$ for which $\phi$ is not an immersion is a $K$-stable analytic
subset of $\red(U)$ which does not intersect $\pi_{\rm
r}^{-1}(y)$. Hence we may shrink $U$ and $Z$ such that they are still
Stein and $\phi\vert U\colon U \to Z$ is a closed embedding. \qed

\nbigskip
In the sequel, let $Z$ be a reduced Stein $K$-space and assume that
$X$ is a $K$-subspace of $Z$. Then $X$ is defined by a continuous
coherent $K$-sheaf ${\cal I}$ of ideals in ${\cal O}_Z$. Let $\kappa
\colon Z \to Z \qq K$ denote the semistable quotient.

\nbigskip
{\bf Proposition 2.} {\it The  ideal $(\kappa_*{\cal I})^K$
defines a closed subspace $X \qq K$ of $Z \qq K$ and the restriction
$\pi := \kappa \vert X \colon X \to X \qq K$ is a semistable quotient
for $X$. Moreover, $\red(\pi)\colon \red(X) \to \red(X \qq K)$ is a
semistable quotient for $\red(X)$.}

\nbigskip
{\it Proof.} By Proposition 1, $(\kappa_*{\cal I})^K$ is coherent and
hence $X\qq K$ is a complex subspace of $Z \qq K$. The fact that $\pi$ and
$\red(\pi)$ are semistable quotients now follows from applying
Lemma~2 to the following two exact sequences of $K$-sheaves:
$$ 0 \to {\cal I} \to {\cal O}_Z \to {\cal O}_X \to 0, \qquad %
0 \to \sqrt{\cal I} \to {\cal O}_Z \to {\cal O}_{\red X} \to
0. \eqno{{\scriptstyle\square}} $$

\nbigskip
Now assume that the ideal ${\cal I}$ has an extension to $Z^c$ by a
holomorphic $K^\C$-ideal ${\cal I}^c$ of ${\cal O}_{Z^c}$. Note that
for the $K^\C$-subspace $X^c$ of $Z$ defined by ${\cal I}^c$ we have
$X^c = K^\C \mal X$. In particular, [H], Section 3.3, implies that
$\red(X^c)$ is the complexification of $\red(X)$. This statement holds
also for the corresponding non-reduced spaces:

\nbigskip
{\bf Proposition 3.} {\it Let $Q \subset X \qq K$ be open and 
$U := \pi^{-1}(Q)$. Then the open subspace $U^c := K^\C
\mal U$ of $X^c$ is the complexification of $U$.}

\nbigskip
{\it Proof.}
We have to show that $U^c$ is universal with respect to
$K$-equivariant holomorphic maps $\phi \colon U \to Y$ into
holomorphic $K^\C$-spaces $Y$.

\smallskip
As a set $U$ is of the form $U = U^c \cap X$. In particular, $U$
is orbit-convex (see [H], 3.2) and hence $\red(U^c)$ is the complexification of
$\red(U)$. So the map $\phi$ extends to a $K^\C$-equivariant
continuous map $\phi^c \colon U^c\to Y$.

\smallskip
Let $\phi^0 \colon {\cal O}_Y \to \phi_*{\cal O}_U$ denote the
comorphism of $\phi$. Then we define the comorphism $(\phi^c)^0
\colon {\cal O}_Y \to \phi^c_*{\cal O}_{U^c}$ as follows: For $y \in
Y$ and $x \in U$ with $\phi^c(x) = y$ choose a $g_0 \in K^\C$
with $g_0 \mal x \in U$ and set

$$(\phi^c)^0({\bf f}) := g_0^{-1} \mal \phi^0(g_0 \mal {\bf f})$$
for every germ ${\bf f} \in {\cal O}_Y$ at $y$. We have to show that this
definition does not depend on the choice of $g_0$. Set
$$ N(x) := \{ g \in K^\C; \; g \mal x \in U \}.$$
Recall that $U$ is an orbit convex subset of $U^c$ and therefore
$N(x)$ is connected mod $K$ (see [H], 1.5, 3.2 and 6.6). Now,
let ${\bf h} := g_0 \mal {\bf f}$ and consider the set
$$M := \{g \in N(x); \; %
(gg_0^{-1})^{-1} \mal \phi^0(gg_0^{-1}\mal {\bf h}) %
= \phi^0({\bf h})\} . $$
Then $Kg_0 \subset M$ holds. Moreover, $M$ is closed in $N(x)$ because
the $K^\C$-sheaves ${\cal O}_{X^c}$ and ${\cal O}_Y$ are
continuous. We claim that $M$ is also open in $N(x)$. Using the
holomorphy of the $K^\C$-sheaves ${\cal O}_{X^c}$ and ${\cal O}_Y$
this is seen as follows: If $g_1\in M$ and $g_2:=g_1g_0^{-1}$, then
$$\phi^0({\bf h})=(gg_2)^{-1}\mal \phi^0(gg_2\mal {\bf h}) 
\leqno (2)$$
holds for all $g\in K$. After representing ${\bf h}$ by some section
$h$ defined in a neighborhood of $\phi(g_0\mal x)$ and restricting
$(2)$ to a sufficiently small neighborhood of $g_0\mal x$, the right
hand side of $(2)$ depends holomorphically on $g$.

\smallskip
Now the Identity Theorem yields that $(2)$ is satisfied for $g$
in some neighborhood of $e \in K^\C$. This implies that $M$ is open in
$N(x)$. Consequently $M = N(x)$ holds, i.e., $(\phi^c)^0$ is well
defined.
\qed

\nbbigskip\goodbreak\noindent
{\bf 5. Proof of the Theorems}

\nbigskip
{\it Proof of the Quotient Theorem.} Assume first that a semistable
quotient $\pi \colon X \to X \qq K$ exists. Consider the exact sequence
$0 \to {\cal N}_X \to {\cal O}_X \to {\cal O}_X / {\cal N}_X \to 0 $
of $K$-sheaves. By Lemma 2, the associated sequence of sheaves of
invariants on $X \qq K$ is also exact. This implies that
$\red(\pi)\colon \red(X) \to \red(X \qq K)$ is a semistable quotient
for $\red(X)$.

\smallskip
Now assume that there is a semistable quotient $\pir \colon
\red(X) \to \red(X) \qq K$. On the level of topological spaces set $X
\qq K := \red(X) \qq K$ and endow $X \qq K$ with the structure sheaf
${\cal O}_{X \qq K} := {\pir}_*({\cal O}_X)^K$. Then it follows
from Proposition 2 and the Local Embedding Lemma that $X \qq K$ is a
complex space. Moreover, by construction, $\pi := \pir$,
interpreted as a morphism of the ringed spaces $X$ and $X \qq
K$, is the quotient map.
\qed

\nbigskip
{\it  Proof of the Coherence Theorem.} By the Local Embedding Lemma and
Proposition~2, we may assume that $X$ is a $K$-stable closed subspace
of a reduced Stein $K$-space $Z$. Noting that the trivial extension of
${\cal S}$ to $Z$ is a continuous coherent $K$-sheaf, we obtain the
assertion from Proposition~1.
\qed

\nbigskip
{\it Proof of the Complexification Theorem.} Let $\pi : X \to X \qq K$
be the semistable quotient for $X$. Choose a cover of $X \qq K$ by
Stein open sets $Q_i$ as in the Local Embedding Lemma. According to
the Local Extension Lemma we may assume that every $U_i :=
\pi^{-1}(Q_i)$ satisfies the assumptions of Proposition~3.

\smallskip
Let $U_i^c$ be the complexification of $U_i := \pi^{-1}(Q_i)$. Note
that every complexification $(U_i \cap U_j)^c$ of $U_i \cap U_j$ is
contained as an open subspace in both, $U_i^c$ and
$U_j^c$. Consequently the $U_i^c$ can be glued together over $X \qq K$
to a complexification $X^c$ of $X$.

\smallskip
By construction, $\red(X^c)$ is the complexification of $\red(X)$.
So we obtain $\red(\pi^c) = (\red \pi)^c$ for the $K$-invariant
holomorphic $\pi^c \colon X^c \to X \qq K$ extension of $\pi$. In
particular, $\pi^c$ is a Stein map and hence a semistable quotient.
\qed

\nbigskip
Note that, by construction, the complexification $X^c$ satisfies
the technical assumptions i) and ii), made in Section 3.

\nbigskip
{\it  Proof of the Extension Theorem.} By the Homomorphism
Lemma, we only have to prove the existence of a
$K^\C$-extension. According to the Local Extension Lemma we can
cover $X \qq K$ by open sets $Q_i$ such that on each $U_i :=
\pi^{-1}(Q_i)$ the restriction ${\cal S}_i := {\cal S}\vert U_i$ has
an $K^\C$-extension ${\cal S}_i^c$ to $U_i^c = (\pi^c)^{-1}(Q_i)$.

\smallskip
The Homomorphism Lemma yields glueing homomorphisms ${\cal S}_i^c
\vert U_i^c \cap U_j^c \to {\cal S}_j^c \vert U_i^c \cap U_j^c $ that
extend the identity map ${\cal S}_i \vert U_i \cap U_j \to {\cal S}_j
\vert U_i \cap U_j$. Hence the ${\cal S}_i^c$ can be glued together
to a $K^\C$-sheaf ${\cal S}^c$ on $X^c$. It is straightforward to
check that ${\cal S}^c$ is a $K^\C$-extension of ${\cal S}$.
\qed

\nbbigskip\goodbreak\noindent
{\bf  References}
\nobreak
\nbigskip\nobreak
\beginlist{1.9cm}
\litem{[G-R]} Grauert, H.; Remmert, R.: Theorie der Steinschen
               R\"aume. 
              Heidelberg: Sprin\-ger 1977
\litem{[HC]} Harish-Chandra: {\it Discrete Series for semisimple Lie
              groups II.} 
             Acta. Math. 116, 1-111 (1966)
\litem{[H]} Heinzner, P.: {\it Geometric Invariant Theory on Stein
              Spaces.} 
            Math. Ann. 289, 631-662 (1991)
\litem{[H-M-P]} Heinzner, P.; Migliorini, L.; Polito, M.:
               {\it Semistable Quotients.} 
               To appear, Annali della Scuola Normale Superiore di
              Pisa (1997)
\litem{[K]} Kaup, W.: {\it Infinitesimale Transformationsgruppen auf
              komplexen R\"{a}umen.}
            Math. Ann. 160, 72-92 (1965)
\litem{[R]} Roberts, M.: {\it A Note on Coherent $G$-Sheaves.} 
            Math. Ann. 275, 573-582 (1986)
\litem{[Sch]} Schwarz, G.: {\it Lifting smooth homotopies of orbit
              spaces.}
              Publ. Math. IHES 51, 37-135  (1980) 
\litem{[S]} Siu, Y.-T.: {\it Every Stein subvariety admits a Stein
              neighborhood.} Inv. Math. 38, 89-100 (1976)

\endlist

\nbbigskip
\settabs 2 \columns

\+ J\"urgen Hausen & Peter Heinzner \cr
\+ Fakult\"at f\"ur Mathematik und & Department of Mathematics \cr
\+ Informatik,  Universit\"at Konstanz & Brandeis University \cr 
\+ 78457 Konstanz & Waltham, MA 02254-9110 \cr
\+ Germany & USA \cr
\+ e-mail: Juergen.Hausen@uni-konstanz.de & e-mail:
heinzner@max.math.brandeis.edu \cr

\end